\documentstyle[11pt,fullpage,twoside,amssym,saks]{article}

\def\XXint#1#2#3{{\setbox0=\hbox{$#1{#2#3}{\int}$}
     \vcenter{\hbox{$#2#3$}}\kern-.5\wd0}}

\begin{document}
\markboth{\centerline{B. Bojarski, V. Gutlyanski and V. Ryazanov}}
{\centerline{ON INTEGRAL CONDITIONS FOR THE GENERAL BELTRAMI
EQUATIONS}}
\def\kohta #1 #2\par{\par\noindent\rlap{#1)}\hskip30pt
\hangindent30pt #2\par}
\def\esssup{\operatornamewithlimits{ess\,sup}}
\def\tomes{\mathop{\longrightarrow}\limits^{mes}}
\def\ts{\textstyle}
\def\I{\roman{Im}}
\def\mes{\mbox{\rm mes}}
\def\Rm{{{\Bbb R}^m}}
\def\Rn{{{\Bbb R}^n}}
\def\Rk{{{\Bbb R}^k}}
\def\R3{{{\Bbb R}^3}}
\def\lR{{\overline {{\Bbb R}}}}
\def\lC{{\overline {{\Bbb C}}}}
\def\lQ{{\overline {{Q}}}}
\def\lz{{\overline {{z}}}}
\def\lD{{\overline {{D}}}}
\def\lR{{\overline {{\Bbb R}}}}
\def\lRn{{\overline {{\Bbb R}^n}}}
\def\lRm{{\overline {{\Bbb R}^m}}}
\def\lBn{{\overline {{\Bbb B}^n}}}
\def\Bn{{{\Bbb B}^n}}
\def\R{{\Bbb R}}
\def\N{{\Bbb N}}
\def\Z{{\Bbb Z}}
\def\C{{\Bbb C}}
\def\B{{\Bbb B}}
\def\Di{{\Bbb D}}
\def\e{{\varepsilon}}
\def\L{{\Lambda}}
\def\l{{\lambda}}
\def\f{{\varphi}}
\def\F{{\Phi}}
\def\x{{\chi}}
\def\d{{\delta }}
\def\D{{\Delta }}
\def\c{{\circ }}
\def\tg{{\tilde{\gamma}}}
\def\a{{\alpha }}
\def\b{{\beta }}
\def\p{{\psi }}
\def\m{{\mu }}
\def\n{{\nu }}
\def\r{{\rho }}
\def\t{{\tau }}
\def\S{{\Sigma }}
\def\O{{\Omega }}
\def\o{{\omega }}
\def\s{{\sigma }}
\def\v{{\vartheta}}
\def\z{{\zeta }}
\def\L{{{Log}}}
\def\E{{{Exp}}}
\def\g{{\gamma }}
\def\G{{\Gamma }}
\def\D{{\Delta }}
\let\text=\mbox
\let\Cal=\cal

\def\cc{\setcounter{equation}{0}
\setcounter{figure}{0}\setcounter{table}{0}}

\overfullrule=0pt

\def\eqb{\begin{equation}}
\def\eqe{\end{equation}}
\def\eb{\begin{eqnarray}}
\def\ee{\end{eqnarray}}
\def\ebnn{\begin{eqnarray*}}
\def\eenn{\end{eqnarray*}}
\def\db{\begin{displaystyle}}
\def\de{\end{displaystyle}}
\def\tb{\begin{textstyle}}
\def\te{\end{textstyle}}
\def\exb{\begin{ex}}
\def\exe{\end{ex}}
\def\bth{\begin{theo}}
\def\eth{\end{theo}}
\def\bcor{\begin{corol}}
\def\ecor{\end{corol}}
\def\blem{\begin{lemma}}
\def\elem{\end{lemma}}
\def\brem{\begin{rem}}
\def\erem{\end{rem}}
\def\bpr{\begin{propo}}
\def\epr{\end{propo}}
\title{{\bf ON  INTEGRAL CONDITIONS FOR THE\\ GENERAL BELTRAMI EQUATIONS}}

\author{{\bf B. Bojarski, V. Gutlyanski and V. Ryazanov}}
\maketitle

\large \abstract Under integral restrictions on dilatations, it is
proved existence theorems for the degenerate Beltrami equations with
two characteristics ${\overline {\partial}}f\, =\, \mu {{\partial f}
} + \nu {\overline {\partial f} }$ and, in particular, to the
Beltrami equations of the second type ${\overline {\partial}}f\, =\,
\nu {\overline {\partial f} }$ that play a great role in many
problems of mathematical physics and to the so--called reduced
Beltrami equations ${\overline {\partial}}f\, =\, \l\ {\rm Re}\,
{{\partial f}}$ that also have significant applications.
\endabstract


\large
\cc
\section{Introduction} The existence problem for the Beltrami equations with two
characteristics \eqb \label{eq1.3} f_{\overline{z}}\, =\, \mu
(z)\cdot{f_z} + \nu (z)\cdot \overline {f_z} \eqe where  $|\mu (z)|+
|\nu (z)| < 1$ a.e. was solved first in the case of the bounded
dilatations \eqb \label{eq1.4} K_{\mu , \nu}(z)\ \colon =\
\frac{1+|\mu (z)|+|\nu (z)|}{1-|\mu (z)|-|\nu (z)|}\eqe in
\cite{Bo$_1$}, Theorem 5.1.\bigskip

Recently in \cite{BGR$_1$}, the existence of homeomophic solutions
in $W^{1,s}_{loc}$ for all $s\in [1,2)$ to the equation
(\ref{eq1.3}) was stated in the case when $K_{\mu , \nu}$ had a
majorant $Q$ in the class BMO, bounded mean oscillation by
John--Nirenberg, see \cite{JN}. Note that $L^{\infty}\subset\ $BMO$\
\subset L^p_{loc}$ for all $p\in [1,\infty)$. Our last paper
\cite{BGR$_2$} was devoted to the study of more general cases when
$K_{\mu , \nu}\in L^{1}_{loc}$ and when $K_{\mu , \nu}$ had a
majorant $Q$ in the class FMO, finite mean oscillation by
Ignat'ev--Ryazanov, see \cite{IR}. Moreover, the paper
\cite{BGR$_2$} contained one new criterion of the Lehto type that is
the base for the further development of the theory of the degenerate
Beltrami equations (\ref{eq1.3}) with integral constraints on the
dilatation $K_{\mu , \nu}$ in the present paper, see the next
section.
\medskip

In the theory of quasiconformal mappings, it is well-known the role
of the Beltrami equations of the first type \eqb \label{eq1.1}
f_{\overline{z}}\, =\, \mu (z)\cdot f_z \eqe where $
f_{\overline{z}} = {\overline {\partial}}f = (f_x+if_y)/2 ,$ $
f_{{z}} = \partial f = (f_x-if_y)/2, $ $ z = x+iy, $ and $f_x$ and
$f_y$ are partial derivatives of $f=u+iv$ in the variables $x$ and
$y,$ respectively, and $\m: D\to \C$ is a measurable function with
$|\m (z)| < 1$ a.e.,  see e.g. \cite{Ah$_1$}, \cite{Bel},
\cite{Bo$_1$} and \cite{LV} where the existence problem was resolved
for the uniformly elliptic case when $\Vert\mu\Vert_{\infty}<1$. The
existence problem for degenerate Beltrami equations (\ref{eq1.1})
with unbounded dilatations \eqb \label{eq1.4a} K_{\mu}(z)\ \colon =\
\frac{1+|\mu (z)|}{1-|\mu (z)|}\eqe is currently an active area of
research, see e.g. \cite{AIM},
 \cite{BJ$_1$}--\cite{BJ$_2$}, \cite{Ch}, \cite{Da},
\cite{GMSV}, \cite{IM}, \cite{Le}, \cite{MM}, \cite{MMV},
\cite{MRSY}, \cite{MS}, \cite{RSY$_1$}--\cite{RSY$_4$}, \cite{SY},
\cite{Tu}, and \cite{Ya}. The study of such homeomorphisms was
started in the frames of the theory of the so--called mean
quasiconformal mappings, see e.g. \cite{Ah$_2$}, \cite{Bil},
\cite{Go}, \cite{GK}, \cite{GR}, \cite{Kr}, \cite{Kru}, \cite{Kud},
\cite{Ku}, \cite{Per}, \cite{Pes}, \cite{Rya}, \cite{Str} and
\cite{UV}.\medskip

On the other hand, the Beltrami equations of the second type \eqb
\label{eq1.2} f_{\overline{z}}\, =\, \nu (z)\cdot \overline {f_z}
\eqe play a great role in many problems of mathematical physics, see
e.g. \cite{KK}. Hence the research of the equations (\ref{eq1.3}) is
so actual. \medskip

Recall that a function $f: D\to\C$ is {\bf absolutely continuous
on lines}, abbr. $f\in${\bf ACL}, if, for every closed rectangle
$R$ in $D$ whose sides are parallel to the coordinate axes, $f|R$
is absolutely continuous on almost all line segments in $R$ which
are parallel to the sides of $R.$ In particular, $f$ is ACL
(possibly modified on a set of Lebesgue measure zero) if it
belongs to the Sobolev class $W^{1,1}_{loc}$ of locally integrable
functions with locally integrable first generalized derivatives
and, conversely, if $f\in $ ACL has locally integrable first
partial derivatives, then $f\in W^{1,1}_{loc}$, see e.g. 1.2.4 in
\cite{MP}.\medskip

If $f: D\to \C$ is a homeomorphic ACL solution of the Beltrami
equation (\ref{eq1.3}) with $K_{\mu , \nu}\in L^{1}_{loc}(D)$, then
$f\in W^{1,1}_{loc}(D)$, furthermore, if $K_{\mu , \nu}\in
L^{p}_{loc}(D)$, $p\in [1,\infty]$, then $f\in W^{1,s}_{loc}(D)$
where $s=2p/(p+1)$. Indeed, if $f\in$ ACL, then $f$ has partial
derivatives $f_x$ and $f_y$ a.e. and, for a sense-preserving ACL
homeomorphism $f: D\to\C ,$ the Jacobian $J_f(z) =
|f_z|^2-|f_{\lz}|^2$ is nonnegative a.e. and, moreover, \eqb
\label{eq1.6} |\ \overline{\partial}f|\ \leq\ |\ \partial f|\ \leq\
|\
\partial f|\ +\ |\ \overline{\partial}f|\ \leq\ Q^{1/2}(z)\cdot
J_f^{1/2}(z)\ \ \ \  \ \ \  a.e.\eqe Recall that if a
homeomorphism $f:D\to\C$ has finite partial derivatives a.e., then
\eqb \label{eq1.7} \int\limits_B\ J_f(z)\ dxdy\ \le\ |f(B)| \eqe
for every Borel set $B\subseteq D$, see e.g. Lemma III.3.3 in
\cite{LV}. Consequently, applying successively the H\"older
inequality and the inequality (\ref{eq1.7}) to (\ref{eq1.6}), we
get that \eqb \label{eq1.8} \Vert\partial f\Vert_s\ \le \ \Vert
K_{\mu , \nu}\Vert^{1/2}_p \cdot |f(C)|^{1/2} \eqe where
$\Vert\cdot\Vert_s$ and $\Vert\cdot\Vert_p$ denote the $L^s-$ and
$L^p-$norm in a compact set $C\subset D$, respectively.
\medskip

In the classical case when $\Vert \m\Vert_{\infty} < 1,$
equivalently, when $K_{\m}\in L^{\infty}(D),$ every ACL homeomorphic
solution $f$ of the Beltrami equation (\ref{eq1.1}) is in the class
$W^{1,2}_{loc}(D)$ with $f^{-1}\in W^{1,2}_{loc}(f(D))$. In the case
$\Vert \m\Vert_{\infty} = 1$ with $K_{\m}\le Q\in \mbox{BMO},$ again
$f^{-1}\in W^{1,2}_{loc}(f(D))$ and $f$ belongs to
$W^{1,s}_{loc}(D)$ for all $1\le s <2$ but already not necessarily
to $W^{1,2}_{loc}(D)$. However, there is a number of degenerate
Beltrami equations (\ref{eq1.1}) for which there exist homeomorphic
solutions $f$ of the class $W^{1,1}_{loc}(D)$ with $f^{-1}\in
W^{1,2}_{loc}(f(D))$.
\medskip

Following \cite{BGR$_2$}, we call a homeomorphism $f\in
W^{1,1}_{loc}(D)$ a {\bf regular solution} of (\ref{eq1.3}) if $f$
satisfies (\ref{eq1.3}) and $J_f(z)\ne 0$ a.e. Note that by
\cite{HK} $f^{-1}\in W^{1,2}_{loc}(f(D))$ for such solutions if
$K_{\mu , \nu}\in L^1_{loc}(D)$.
\medskip

\cc
\section{Preliminaries}   The following theorem was recently
established in the work \cite{BGR$_2$}.

\bth{} \label{th3.2C}  Let  $D$  be a domain in $\C$
 and let $\mu$
and $\nu : D\to\C$ be measurable functions with $|\mu (z)|+|\nu (z)|
< 1$ a.e. and $K_{\mu , \nu}\in L^1_{loc}(D).$ Suppose that
 \eqb \label{eq4.30C}
 \int\limits_{0}^{\d(z_0)}\frac{dr}{rk_{z_0}( r)}\ =\ \infty
 \ \ \  \ \ \ \ \ \ \forall\ z_0\in D  \eqe
where $\d(z_0)<dist\, (z_0,\partial D)$ and $k_{z_0}(r)$ is the
average of $K_{\mu , \nu}(z)$ over the circle $|z-z_0|=r.$ Then the
Beltrami equation (\ref{eq1.3}) has a regular solution. \eth

In general, in the Beltrami equation theory in the plane as well as
 in the theory of space mappings, the integral conditions of the Lehto type
\begin{equation}\label{eq1} \int\limits_{0}^{1}\ \frac{dr}{rq(r)}\
=\ \infty \end{equation} are often met where the function $Q$ is
given say in the unit ball ${\Bbb B}^n=\{ x\in{\Bbb R}^n: |x|<1\}$
and $q(r)$ is the average of the function $Q(z)$ over the sphere
$|x|=r$, see e.g. \cite{AIM}, \cite{BGR$_2$}, \cite{Ch},
\cite{GMSV}, \cite{Le}, \cite{MRSY}, \cite{MS}, \cite{Per},
\cite{RSY$_1$}--\cite{RSY$_4$}, \cite{Zo$_1$} and
\cite{Zo$_2$}.\bigskip

On the other hand, in the theory of  mappings called quasiconformal
in the mean, conditions of the type
\begin{equation}\label{eq2} \int\limits_{{\Bbb B}^n} \Phi
(Q(x))\ dx\  <\ \infty\end{equation} are standard for various
characteristics $Q$ of these mappings, see e.g. \cite{Ah$_2$},
\cite{Bil}, \cite{Go}, \cite{GR}, \cite{Kr}--\cite{Ku}, \cite{Per},
\cite{Rya} and \cite{Str}.\bigskip

In this connection, in the paper \cite{RSY$_4$} it was established
in\-ter\-con\-nec\-tions between a series of integral conditions on
the function $\Phi$ and between (\ref{eq1}) and (\ref{eq2}), cf.
also \cite{BJ$_2$} and \cite{GMSV}. We give here these conditions
for $\Phi$ under which (\ref{eq2}) implies (\ref{eq1}).\medskip

Further we use the following notion of the inverse function for
mo\-no\-to\-ne functions. For every non-decreasing function
$\Phi:[0,\infty ]\to [0,\infty ] ,$ the {\bf inverse function}
$\Phi^{-1}:[0,\infty ]\to [0,\infty ]$ can be well defined by
setting
\begin{equation}\label{eq5.5CC} \Phi^{-1}(\tau)\ =\
\inf\limits_{\Phi(t)\ge \tau}\ t\ .
\end{equation} As usual, here $\inf$ is equal to $\infty$ if the set of
$t\in[0,\infty ]$ such that $\Phi(t)\ge \tau$ is empty. Note that
the function $\Phi^{-1}$ is non-decreasing, too.

\brem\label{rmk3.333} It is evident immediately by the definition
that \eqb\label{eq5.5CCC} \F^{-1}(\F(t))\ \le\ t\ \ \ \ \ \ \ \
\forall\ t\in[ 0,\infty ] \eqe with the equality in (\ref{eq5.5CCC})
except intervals of constancy of the function $\Phi(t)$. \erem

Further, the integral in (\ref{eq333F}) is un\-der\-stood as the
Lebesgue--Stieltjes integral and the integrals in (\ref{eq333Y}) and
(\ref{eq333B})--(\ref{eq333A}) as the ordinary Lebesgue integrals.
In (\ref{eq333Y}) and (\ref{eq333F}) we complete the definition of
integrals by $\infty$ if $\Phi(t)=\infty ,$ correspondingly,
$H(t)=\infty ,$ for all $t\ge T\in[0,\infty) .$

\bth{} \label{pr4.1aB} Let $\F:[0,\infty ]\to [0,\infty ]$ be a
non-decreasing function and set \eqb\label{eq333E} H(t)\ =\ \log
\F(t)\ .\eqe

Then the equality \eqb\label{eq333Y} \int\limits_{\D}^{\infty}
H'(t)\ \frac{dt}{t}\ =\ \infty   \eqe implies the equality
\eqb\label{eq333F} \int\limits_{\D}^{\infty} \frac{dH(t)}{t}\ =\
\infty  \eqe and (\ref{eq333F}) is equivalent to \eqb\label{eq333B}
\int\limits_{\D}^{\infty}H(t)\ \frac{dt}{t^2}\ =\ \infty \eqe for
some $\D>0,$ and (\ref{eq333B}) is equivalent to every of the
equalities: \eqb\label{eq333C}
\int\limits_{0}^{\d}H\left(\frac{1}{t}\right)\ {dt}\ =\ \infty \eqe
for some $\d>0,$ \eqb\label{eq333D} \int\limits_{\D_*}^{\infty}
\frac{d\eta}{H^{-1}(\eta)}\ =\ \infty \eqe for some $\D_*>H(+0),$
\eqb\label{eq333A} \int\limits_{\d_*}^{\infty}\ \frac{d\t}{\t
\F^{-1}(\t )}\ =\ \infty \eqe for some $\d_*>\F(+0).$
\medskip

Moreover, (\ref{eq333Y}) is equivalent  to (\ref{eq333F}) and hence
(\ref{eq333Y})--(\ref{eq333A})
 are equivalent each to other  if $\F$ is in addition absolutely continuous.
In particular, all the conditions (\ref{eq333Y})--(\ref{eq333A}) are
equivalent if $\F$ is convex and non--decreasing. \eth

It is necessary here to give one more explanation. From the right
hand sides in the conditions (\ref{eq333Y})--(\ref{eq333A}) we have
in mind $+\infty$. If $\Phi(t)=0$ for $t\in[0,t_*]$, then
$H(t)=-\infty$ for $t\in[0,t_*]$ and we complete the definition
$H'(t)=0$ for $t\in[0,t_*]$. Note, the conditions (\ref{eq333F}) and
(\ref{eq333B}) exclude that $t_*$ belongs to the interval of
integrability because in the contrary case the left hand sides in
(\ref{eq333F}) and (\ref{eq333B}) are either equal to $-\infty$ or
indeterminate. Hence we may assume in (\ref{eq333Y})--(\ref{eq333C})
that $\Delta>t_0$ where $t_0\colon =\sup\limits_{\Phi(t)=0}t$,
$t_0=0$ if $\Phi(0)>0$, and $\delta<1/t_0$, correspondingly.\bigskip

Finally, we give the connection of the above conditions with the
condition of the Lehto type (\ref{eq4.30C}).
\medskip

Recall that a function  $\psi :[0,\infty ]\to [0,\infty ]$ is called
{\bf convex} if $\psi (\lambda t_1 + (1-\lambda) t_2)\le\lambda\psi
(t_1)+ (1-\lambda)\psi (t_2)$ for all $t_1$ and $t_2\in[0,\infty ]$
and $\lambda\in [0,1]$.\medskip

In what follows, $\Di$ denotes the unit disk in the complex plane
$\C$, \eqb\label{eq5.5Cf} \Di\ =\ \{\ z\in\C:\ |z|\ <\ 1\ \}\ .\eqe

\bth{} \label{th5.555} Let $Q:\Di\to [0,\infty ]$ be a measurable
function such that \eqb\label{eq5.555} \int\limits_{\Di} \F (Q(z))\
dxdy\  <\ \infty\eqe where $\F:[0,\infty ]\to [0,\infty ]$ is a
non-decreasing convex function such that \eqb\label{eq3.333a}
\int\limits_{\d_0}^{\infty}\ \frac{d\t}{\t \F^{-1}(\t )}\ =\ \infty
\eqe for some $\d_0\ >\ \t_0\ \colon =\ \F(0).$ Then
\eqb\label{eq3.333A} \int\limits_{0}^{1}\ \frac{dr}{rq(r)}\ =\
\infty \eqe where $q(r)$ is the average of the function $Q(z)$ over
the circle $|z|=r$. \eth

Finally, combining Theorems \ref{pr4.1aB} and \ref{th5.555} we
obtain the following conclusion.

\bcor \label{cor555} If $\F:[0,\infty ]\to [0,\infty ]$ is a
non-decreasing convex function and $Q$ satisfies the condition
(\ref{eq5.555}), then every of the conditions
(\ref{eq333Y})--(\ref{eq333A}) implies (\ref{eq3.333A}). \ecor

\cc
\section{Existence theorems} Immediately on the base of Theorem \ref{th3.2C} and
Corollary \ref{cor555}, we obtain the next significant result.

\bth{} \label{th4.111a} Let  $D$  be a domain in $\C$ and let $\mu$
and $\nu : D\to\C$ be measurable functions with $|\mu (z)|+|\nu (z)|
< 1$ a.e. such that \eqb\label{eq4.2} \int\limits_{{D}} \Phi
(K_{\m,\n}(z))\ dxdy\  <\ \infty \eqe where $\F:[0,\infty ]\to
[0,\infty ]$ is a non-decreasing convex function. If $\Phi$
satisfies at least one of the conditions
(\ref{eq333Y})--(\ref{eq333A}), then the Beltrami equation
(\ref{eq1.3}) has a regular solution. \eth

\brem\label{rmk5.1} The condition (\ref{eq4.2}) can be also
localized to neighborhoods $U_{z_0}$ of points $z_0\in D$ with
$\Phi=\Phi_{z_0}$ under the same conditions on the functions
$\Phi_{z_0}$. If $\infty\in D$, then the condition (\ref{eq4.2}) for
$K_{\m,\n}(z)$ at $\infty\in D$ should be understood as the
corresponding condition for $K_{\m,\n}({1}/{\overline z})$ at $0.$
The latter condition can also be rewritten explicitly in terms of
$K_{\m,\n}(z)$ itself after the inverse change of variables
$z\longmapsto{1}/{\overline z}$ in the form \eqb\label{eq3.4R}
\int\limits_{{U_{\infty}}} \Phi_{\infty} (K_{\m,\n}(z))\
\frac{dxdy}{|z|^4}\ <\ \infty\ . \eqe If the domain $D$ is
unbounded, then it is better to use the global condition
\eqb\label{eq3.4R} \int\limits_{{D}} \Phi (K_{\m,\n}(z))\
\frac{dxdy}{(1+|z|^2)^2}\  <\ \infty \eqe instead of the condition
(\ref{eq4.2}). The latter means the integration of the function
$\Phi\circ K_{\m,\n}$ in the spherical area.
\medskip

We may assume in the above theorem that the functions $\F_{z_0}(t)$
and $\F(t)$ are not convex and non--decreasing on the whole segment
$[0,\infty]$ but only on a segment $[T,\infty]$ for some
$T\in(1,\infty)$. Indeed, every function
$\F:[0,\infty]\to[0,\infty]$ which is convex and non-decreasing on a
segment $[T,\infty]$, $T\in(0,\infty)$, can be replaced by a
non-decreasing convex function $\F_T:[0,\infty]\to[0,\infty]$ in the
following way. We set $\F_T(t)\equiv 0$ for all $t\in [0,T]$,
$\F(t)=\f(t)$, $t\in[T,T_*]$, and $\F_T\equiv \F(t)$,
$t\in[T_*,\infty]$, where $\t=\f(t)$ is the line passing through the
point $(0,T)$ and touching upon the graph of the function $\t=\F(t)$
at a point $(T_*,\F(T_*))$, $T_*\ge T$. For such a function we have
by the construction that $\F_T(t)\le \F(t)$ for all $t\in[1,\infty]$
and $\F_T(t)=\F(t)$ for all $t\ge T_*$. \erem

The equation of the form \eqb \label{eq6.1} f_{\overline{z}}\, =\,
\l (z)\ {\rm Re}\, f_z \eqe with $|\l(z)|<1$ a.e. is called the {\bf
reduced Beltrami equation}, see e.g. \cite{AIM},
\cite{Bo$_1$}--\cite{Bo$_3$} and \cite{Vo}. The equation
(\ref{eq6.1}) can be rewritten as the equation (\ref{eq1.3}) with
\eqb \label{eq6.2} \m(z)\ =\ \n (z)\ =\ \frac{\l(z)}{2} \eqe and
then \eqb \label{eq6.3} K_{\mu , \nu}(z)\ =\ K_{\l}(z)\ \colon =\
\frac{1+|\l (z)|}{1-|\l (z)|}\ .\eqe Thus, we obtain from Theorem
\ref{th4.111a} the following consequence for the reduced Beltrami
equations (\ref{eq6.1}).

\bth{} \label{th4.111aR} Let  $D$  be a domain in $\C$ and let $\l$
be a measurable function with $|\l (z)| < 1$ a.e. such that
\eqb\label{eq2.8aR} \int\limits_{{D}} \Phi (K_{\l}(z))\ dxdy\ <\
\infty \eqe where $\F:[0,\infty ]\to [0,\infty ]$ is a
non-decreasing convex function. If $\Phi$ satisfies at least one of
the conditions (\ref{eq333Y})--(\ref{eq333A}), then the reduced
Beltrami equation (\ref{eq6.1}) has a regular solution. \eth

\brem\label{rmk6.1} Remarks \ref{rmk5.1} are valid for the reduced
Beltrami equation. Moreover, the above results remain true for the
case in (\ref{eq1.3}) when \eqb \label{eq6.14} \n(z)\ =\ \m (z)\
e^{i\theta (z)} \eqe with an arbitrary measurable function
$\theta(z): D\to\R$ and, in particular, for the equations of the
form \eqb \label{eq6.1A} f_{\overline{z}}\, =\, \l (z)\ {\rm Im}\,
f_z \eqe with a measurable coefficient $\l : D\to\C$, $|\l(z)|<1$
a.e., see e.g. \cite{Bo$_1$}--\cite{Bo$_3$}.\medskip

Next, note that Theorem 5.50 from the work \cite{RSY$_4$} for the
Beltrami equations of the first type (\ref{eq1.1}) shows that the
conditions (\ref{eq333Y})--(\ref{eq333A}) are not only sufficient
but also necessary for the general Beltrami equations (\ref{eq1.3})
to have regular solutions.\medskip

Finally, the same is valid for the reduced Beltrami equations
(\ref{eq6.1}) because the examples in the mentioned theorem had the
form $$ f(z)\ =\ \frac{z}{|z|}\ \r(|z|)
$$
where $\r(t)= e^{I(t)}$ and
$$
I(t)\ \colon =\ \int\limits_0^t\ \frac{dr}{rK(r)}\ .
$$
Indeed, setting $z=re^{i\vartheta}$ we have that
$$
\frac{\partial f}{\partial r}\ =\ \frac{\partial f}{\partial z}
\cdot \frac{\partial z}{\partial r}\ +\ \frac{\partial f}{\partial
\overline {z}} \cdot \frac{\partial \overline {z}}{\partial r}\ =\
e^{i\vartheta}\cdot\frac{\partial f}{\partial z} \ +\
e^{-i\vartheta}\cdot\frac{\partial f}{\partial \overline {z}}
$$
and
$$ \frac{\partial f}{\partial \vartheta}\ =\ \frac{\partial
f}{\partial z} \cdot \frac{\partial z}{\partial \vartheta}\ +\
\frac{\partial f}{\partial \overline {z}} \cdot \frac{\partial
\overline {z}}{\partial \vartheta}\ =\
ire^{i\vartheta}\cdot\frac{\partial f}{\partial z} \ -\
ire^{-i\vartheta}\cdot\frac{\partial f}{\partial \overline {z}}
$$
and hence $$ \frac{\partial f}{\partial z}\ =\
\frac{e^{-i\vartheta}}{2}\left(\frac{\partial f}{\partial r} \ +\
\frac{1}{ir}\cdot\frac{\partial f}{\partial \vartheta}\right)\ =\
\frac{1}{2}\left(\frac{\r (r)}{rK(r)} \ +\ \frac{\r(r)}{r}\right)\
=\ \frac{\r(r)}{2r}\ \cdot \frac{1+K(r)}{K(r)}\ >\ 0
$$ and $$ \frac{\partial f}{\partial \overline {z}}\
=\ \frac{e^{i\vartheta}}{2}\left(\frac{\partial f}{\partial r} \ -\
\frac{1}{ir}\cdot\frac{\partial f}{\partial \vartheta}\right)\ =\
\frac{e^{2i\vartheta}}{2}\left(\frac{\r (r)}{rK(r)} \ -\
\frac{\r(r)}{r}\right)\ =\ e^{2i\vartheta}\cdot\frac{\r(r)}{2r}\
\cdot \frac{1-K(r)}{K(r)}\ , $$ i.e.
$$ \l(z)\ =\
e^{2i\vartheta}\cdot\frac{1-K(r)}{1+K(r)}\ =\ -\ \frac{z}{\overline
{z}}\cdot\frac{K(|z|)-1}{K(|z|)+1}
$$
and, consequently, $K_{\l}(z)\ =\ K(|z|)$. \erem

{\bf Acknowledgements.} The research of the third author was
partially supported by Institute of Mathematics of PAN, Warsaw,
Poland, and by Grant F25.1/055 of the State Foundation of
Fundamental Investigations of Ukraine.

\bigskip

\noindent Bogdan Bojarski, Institute of Mathematics of Polish
Academy of Sciences,\\ ul. Sniadeckich 8, P.O. Box 21, 00--956
Warsaw, POLAND\\ Email: {\tt
bojarski@impan.gov.pl}\\

\bigskip

\noindent Vladimir Gutlyanski, Vladimir Ryazanov, Inst. Appl. Math.
Mech., NASU,\\ 74 Roze Luxemburg str., 83114, Donetsk, UKRAINE,\\
Email: {\tt vladimirgut@mail.ru}, {\tt vlryazanov1@rambler.ru}

\end{document}